\documentstyle[11pt]{article}
\newcommand{\A}{{\sl Avoid\,}}
\newcommand{\h}{\hbox{ht}}
\newcommand{\w}{\hbox{wd}}

\newcommand{\qed}{\hfill$\Box$}
\newcommand{\sh}{{\sl shape}}
\newcommand{\av}{{\sl avoid\,}}

\newenvironment{note}[1]{\par\addvspace{\medskipamount}\noindent
                         {\bf {#1}}\sl
                       }{\par\addvspace{\medskipamount}\rm}

\title{Shape Avoiding Permutations}
\author{Ron M. Adin%
\thanks{Department of Mathematics and Computer Science, Bar-Ilan University,
Ramat-Gan 52900, Israel. Email: {\tt radin@math.biu.ac.il}}
$^\ddag$
\and Yuval Roichman%
\thanks{Department of Mathematics and Computer Science, Bar-Ilan University,
Ramat-Gan 52900, Israel. Email: {\tt yuvalr@math.biu.ac.il}}
\thanks{Both authors supported in part by the Israel Science Foundation
and by internal research grants from 
Bar-Ilan University.}}
\date{Submitted: February 20, 2000; Revised: March 21, 2001}

\begin{document}

\maketitle
\begin {abstract}
Permutations avoiding all patterns of a given shape (in the sense of 
Robinson-Schensted-Knuth) are considered.  We show that the shapes of all such 
permutations are contained in a suitable thick hook, and deduce an exponential
growth rate for their number. 
\end{abstract}

\section{Introduction}

\subsection{Outline}

The Robinson-Schensted(-Knuth) correspondence is a bijection between 
permutations in $S_n$ and pairs of standard Young tableaux of the same shape
(and size $n$).  This common shape is called the {\em shape} of the
 permutation.
A permutation $\pi = (\pi_1,\ldots,\pi_n)$ in $S_n$ {\it avoids a permutation}
$\sigma =(\sigma_1,\ldots,\sigma_m)$ in $S_m$ if there is no subsequence 
$(\pi_{i_1},\ldots,\pi_{i_m})$ of $\pi$ such that $\pi_{i_j} > \pi_{i_k}$ iff
$\sigma_j > \sigma_k$ ($\forall j,k$).  $\pi$ {\em avoids a shape} $\mu$ if
it avoids all the permutations of shape $\mu$.

This paper deals with the relation between the property 
``$\pi$ does not avoid a given shape $\mu$'' 
and the property 
``$\lambda=\sh(\pi)$ contains $\mu$ as a sub-shape''.
It turns out that, in general, neither of these properties implies or 
contradicts the other; but in certain important cases, such implications 
do hold.  These cases include, e.g., rectangular shapes and hook shapes
(either for $\lambda$ or for $\mu$). 
These positive results are then applied to get asymptotic bounds related to 
the Stanley-Wilf conjecture on pattern-avoiding permutations (see 
Corollaries 4 and 5 in Subsection 1.2, and Subsection 7.2).  Use is made of
the Berele-Regev asymptotic evaluation of the number of standard Young
tableaux contained in a ``thick hook''.

\medskip

The rest of the paper is organized as follows.
The main results are listed in Subsection 1.2.
Standard notations and necessary background are given in Section 2.
In Section 3 we motivate our investigation by a ``false conjecture".
In Section 4 we show that this ``false conjecture"
is correct for rectangular shapes.
Using this knowledge we consider the general case in Section 5.
Families of shapes, for which an exact
evaluation may be obtained, are presented in Section 6.
Section 7 concludes the paper with final remarks and open problems.

\bigskip

\subsection{Main Results}

For rectangular shapes the following holds.

\begin{note}
\noindent{\bf Theorem 1.}
If $\pi$ is a permutation of rectangular shape $(m^k)$,
and $\mu$ is an arbitrary shape, then:

$\mu$ is the shape of some subsequence of $\pi$ if and only if
 $\mu\subseteq (m^k)$. 
\end{note}

\noindent See Theorem 4.1 below.

\smallskip

Using Theorem 1 we prove the following general result.

\begin{note}
\noindent{\bf Theorem 2.} 
For any permutation $\pi$ in $S_n$
and any partition $\mu=(\mu_1,\dots,\mu_k)$ of $m$ :

If $(\mu_1^k)\subseteq \sh(\pi)$ then $\mu$ is the
 shape of some subsequence of $\pi$. 
\end{note}

\noindent See Theorem 5.1 below.

\smallskip

For hook shapes a stronger result is proved.

\begin{note}
\noindent{\bf Theorem 3.} Let $m$ and $k$ be positive integers
and let $n\ge 4km$.
Then for any hook $\mu=(m,1^{k-1})$ and  any permutation $\pi$ in $S_n$ :

$\pi$ has a subsequence of shape $\mu$ if and only if $\mu\subseteq\sh(\pi)$.
\end{note}

\noindent See Theorem 6.1 below.

\bigskip
 
Denote by $\av_n^\mu$
 the size of the set of all $\mu$-avoiding permutations in $S_n$.  
Combining Theorem 2 with the Berele-Regev asymptotic estimates \cite{BR}
the following bounds are proved.

\begin{note}
\noindent{\bf Corollary 4.} For any fixed partition $\mu=(\mu_1,\dots,\mu_k)$, 
$$
\max\{\h(\mu),\w(\mu)\}\le 
\liminf_{n\rightarrow\infty}{(\av_n^\mu)^{1/2n}}
$$
and
$$
\limsup_{n\rightarrow\infty}{(\av_n^\mu)^{1/2n}} \le
\h(\mu)+\w(\mu) ,
$$
where the height of $\mu$ $\h(\mu):=k-1$, 
and the width of $\mu$ $\w(\mu):=\mu_1-1$.
\end{note}

\noindent See Corollary 5.2 below. 
It should be noted that this result is related to 
the Stanley-Wilf conjecture (see Subsection 7.2).

\smallskip

For hook shapes we have a sharper estimate.

\begin{note}
\noindent{\bf Corollary 5.} For any pair of positive integers
$m$ and $k$
$$
\lim_{n\rightarrow\infty}(\av^{(m,1^{k-1})}_n)^{1/2n}=\max\{m-1,k-1\}.
$$
\end{note}

\noindent See Corollary 6.5 below. 

\bigskip

\section{Preliminaries}

Two classical partial orders on the set of partitions are considered in 
this paper.
Let $\lambda=(\lambda_1,\dots)$ and $\mu=(\mu_1,\dots)$ be two partitions
(not necessarily of the same number). 

We say that $\mu$ is {\it contained} in $\lambda$, denoted 
$\mu\subseteq \lambda$, if
$$
\mu_i\le \lambda_i \qquad (\forall i).
$$

We say that $\mu$ is 
{\it dominated} by  $\lambda$, denoted $\mu\preceq  \lambda$, if 
$$
\sum\limits_{j=1}^i\mu_j\le   \sum\limits_{j=1}^i\lambda_j
\qquad (\forall i) .
$$
Clearly, $\mu\subseteq \lambda\Rightarrow \mu\preceq \lambda$.

\smallskip

\noindent The  partition {\it conjugate} to $\lambda$ is 
$\lambda'=(\lambda'_1,\dots)$,
where 
$\lambda_i'=\max\{j|\lambda_j\ge i\}$;
i.e., the conjugate partition is obtained by interchanging rows and columns in $\lambda$.

\begin{note}
\noindent{\bf Lemma 2.1.} [Md Ch. I (1.11)] If $\lambda$ and $\mu$ are partitions of the same number $n$ then
$$
\mu\preceq \lambda \Leftrightarrow \lambda'\preceq \mu'.
$$
\end{note}

\begin{note}
\noindent{\bf Corollary 2.2.} If $\lambda$ and $\mu$ are partitions of 
the same number $n$, satisfying
$$
\mu\preceq \lambda \hbox{ and } \mu'\preceq \lambda'
$$ 
then $\lambda=\mu$.
\end{note}

\medskip

Define the {\it shape} of a sequence of integers to be the common shape of 
the two tableaux obtained via the Robinson-Schensted-Knuth correspondence.
See [Sa \S 3.3, St \S 7.11]. The following theorem is well known.

\begin{note}
\noindent{\bf Schensted's Theorem.} [Sc]
For any partition $\lambda$ and any
permutation $\pi$ of shape $\lambda$,
the length of the longest increasing subsequence of $\pi$ is
 equal to $\lambda_1$, and
the length of the longest decreasing subsequence of $\pi$ 
is equal to  $\lambda_1'$.
\end{note}

Schensted's Theorem was generalized by Greene.

\begin{note}
\noindent{\bf Greene's Theorem.} [Gr]
Let $\pi$ be a permutation of shape $\lambda=(\lambda_1,\dots,\lambda_t)$.
Then, for all $i$:
$$
\sum\limits_{j=1}^i \lambda_j=
\hbox{ maximal  size of a union of $i$ increasing subsequences in $\pi$},
$$
and
$$
\sum\limits_{j=1}^i \lambda'_j=
\hbox{ maximal size of a union of $i$ decreasing subsequences in $\pi$}.
$$
\end{note}

\section{Motivation}

Let $\mu$ be  a partition of $m$, and let
$C^\mu$ be the set of all permutations
in $S_m$ of shape $\mu$. 
A permutation in $S_n$ is a {\it $\mu$-avoiding permutation} if it
avoids all the permutations in $C^\mu$;
  denote the set of these permutations by $\A_n^\mu$.

The only permutation in $S_m$ having shape $(m)$ is the identity permutation,
i.e., a monotone increasing sequence.
Schensted's Theorem, stated in the previous section, is thus equivalent to 
the following statement.

\begin{note}
\noindent{\bf Fact 3.1.} For any pair of positive integers $m\le n$
$$
\A_n^{(m)}=\bigcup_{\{\lambda\vdash n| (m)\not\subseteq \lambda\}} C^\lambda ,
$$
and similarly for $(1^m)$ instead of $(m)$.
\end{note}

\medskip

In other words, the set of permutations in $S_n$ avoiding $(m)$ is the union 
of all Knuth cells of shapes not containing $(m)$.
One may be tempted to think that this is a general phenomenon.

\begin{note}
\noindent{\bf ``False Conjecture"  (First Version).}
 For any pair of positive integers 
$m\le n$ and any partition $\mu$ of $m$
$$
\A_n^{\mu}=\bigcup_{\{\lambda\vdash n| \mu\not\subseteq \lambda\}} C^\lambda .
$$
\end{note}
 
Equivalently,

\begin{note}
\noindent{\bf ``False Conjecture"  (Second Version).}
For any permutation $\pi\in S_n$ of shape $\lambda$,
the following two assertions hold:
\begin{itemize}
\item[(1)] For any partition $\mu\subseteq \lambda$
there exists a subsequence of $\pi$ of shape $\mu$.
\item[(2)] The shape of any subsequence of $\pi$ is contained in $\lambda$.
\end{itemize}
\end{note}

\noindent Clearly, (1) is equivalent to the inclusion
$$
\A_n^{\mu}\subseteq \bigcup_{\{\lambda\vdash n| \mu\not\subseteq \lambda\}} 
C^\lambda,
$$ 
while (2) is equivalent to the reverse inclusion
$$
 \bigcup_{\{\lambda\vdash n| \mu\not\subseteq \lambda\}} C^\lambda\subseteq 
 \A_n^{\mu} .
$$

Note that Greene's Theorem implies the weaker result that the shape of any
subsequence of $\pi$ is {\bf dominated} by $\lambda$.

Unfortunately, the following examples show that both
 parts of the ``False Conjecture" are false
 in general.

\begin{note}
\noindent{\bf Example 3.2.}
 The permutation $\pi=(65127843)$ has shape $\lambda=(4,2,1^2)$,
 but  has no subsequence of shape $\mu= (4,1^3)$.

\smallskip

\noindent {\bf Example 3.3.}
The permutation
$\pi=(25314)$ has shape $\lambda=(3,1^2)$, but has a subsequence of shape
$\mu= (2^2)$.
\end{note}

Both examples can be extended to shapes $\lambda$ of arbitrarily large size.

\smallskip

A central discovery in this paper is that the above ``False Conjecture" 
is nevertheless correct in some important cases. 
This will be used to deduce asymptotic estimates.

\bigskip

\section{Rectangular Shapes}

 A {\it rectangular shape} is a shape of the form $(m^k)$,
where $m$ and $k$ are positive integers.
In this section we show that the ``False Conjecture'' is true whenever 
$\lambda$ is a rectangular shape.

\begin{note}
\noindent{\bf Theorem 4.1.}
If $\pi$ is a permutation of rectangular shape $(m^k)$,
and $\mu$ is an arbitrary shape, then:

$\mu$ is the shape of some subsequence of $\pi$ if and only if
 $\mu\subseteq (m^k)$. 
\end{note}

\smallskip

In order to prove Theorem 4.1 we need the following consequence of Greene's Theorem.

\begin{note}
\noindent{\bf Lemma 4.2.} Let $\pi$ be a permutation of shape $\lambda$.
\begin{itemize}
\item[(a)] If $\pi$ contains a disjoint union of $k$ increasing subsequences of lengths
$\ell_1\ge \ell_2\ge \cdots\ge \ell_k$ then
$(\ell_1,\dots,\ell_k)\preceq \lambda.$
\item[(b)] If $\pi$ contains a disjoint union of $k$ decreasing subsequences of lengths
$\ell_1\ge \ell_2\ge \cdots\ge \ell_k$ then
$(\ell_1,\dots,\ell_k)\preceq \lambda'.$
\end{itemize}
\end{note}

\smallskip

\noindent{\bf Proof.} By Greene's Theorem, for any $1\le i \le k$
$$
\sum\limits_{j=1}^i \ell_j\le \hbox{maximal size of a union of $i$ increasing subsequences of $\pi$}=\sum\limits_{j=1}^i\lambda_j. 
$$
The proof of the second part is similar.

\qed

\smallskip

The following lemma characterizes permutations having rectangular shape.

\begin{note}
\noindent{\bf Lemma 4.3.}
\begin{itemize}
\item[(a)] A permutation $\pi$ has shape $(m^k)$ if and only if
the following two conditions are simultaneously satisfied:
\begin{itemize}
\item[(a1)] $\pi$ is a disjoint union of $k$ increasing subsequences,
 each of length $m$.
\item[(a2)] $\pi$ is a disjoint union of $m$ decreasing subsequences,
each of length $k$.
\end{itemize}
\item[(b)] If the above conditions hold, then
each of the $k$ increasing subsequences intersects each
of the $m$ decreasing subsequences in  exactly one element.
\end{itemize}
\end{note}

\noindent{\bf Proof.} 

\noindent{\bf (a)} Assume that $\pi$ has shape $\lambda$ and
satisfies conditions (a1) and (a2) of the Lemma. 
By (a1) and Lemma 4.2(a), $(m^k)\preceq \lambda$.
By (a2) and Lemma 4.2(b), $(k^m)\preceq \lambda'$. 
Also $|\lambda|=|(m^k)|=km$, so by
Corollary 2.2, $\lambda=(m^k)$.

\smallskip

In the other direction: By Greene's Theorem,
if $\pi$ has shape $(m^k)$ then it is the 
disjoint union of $k$ increasing subsequences $\alpha_1,\dots,\alpha_k$
of total size $km$. By Schensted's Theorem, each increasing
subsequence of $\pi$ has size at most $m$, and therefore
$|\alpha_1|=\dots=|\alpha_k|=m$.
Similarly, $\pi$ is a disjoint union of $m$ decreasing subsequences
$\beta_1,\dots,\beta_m$ satisfying $|\beta_1|=\dots=|\beta_m|=k$. 

\smallskip

\noindent{\bf (b)} Each increasing subsequence $\alpha_i$ intersects each decreasing subsequence $\beta_j$
in at most one element, and since these $km$ intersections cover all 
elements of $\pi$ they are all nonempty.

\qed

\medskip

\noindent{\bf Proof of Theorem 4.1.}
Let $\pi$ be a sequence of shape $\lambda=(m^k)$.
If $\mu$ is the shape of some subsequence of $\pi$ then this subsequence contains an increasing subsequence of length $\mu_1$. Therefore
 $\mu_1\le\lambda_1=m$. Similarly $\mu_1'\le \lambda_1'=k$, so that
$\mu\subseteq (m^k)$. 

In the other direction:  By Lemma 4.3, $\pi$ is 
a disjoint union
of $k$ increasing subsequences,
 of length $m$ each, say $\alpha_1,\dots,\alpha_k$ (enumerated arbitrarily).
Similarly, $\pi$ is a disjoint union of $m$ decreasing subsequences,
say $\beta_1,\dots,\beta_m$ (of length $k$ each).
Also, each $\alpha_i$ intersects each $\beta_j$ in a unique element; denote it by $P(i,j)$.
Now let $\mu\subseteq (m^k)$, and define 
$\sigma$ to be the subsequence of $\pi$ 
consisting of all elements $P(i,j)$ with $j\le \mu_i$.
We claim that $\sigma$ has shape $\mu$.

Indeed, $\sigma$ intersects $\alpha_i$ in $\mu_i$ elements,
and therefore (by Lemma 4.2(a)) $\mu \preceq \sh(\sigma)$.
Similarly, $\sigma$ intersects $\beta_j$ in $\mu_j'$ elements,
and therefore (by Lemma 4.2(b)) $\mu' \preceq \sh(\sigma)'$.
Since $|\sh(\sigma)| = |\mu|$ by definition, Corollary~2.2 implies that 
$\sh(\sigma)=\mu$ and the proof is complete.

\qed

\medskip

The following theorem is complementary.

\begin{note}
\noindent{\bf Theorem 4.4.} If $\pi$ is a sequence of shape $\lambda$
and $(m^k)\subseteq \lambda$, then there exists a subsequence of 
$\pi$ of shape $(m^k)$.
\end{note}

\noindent In other words: For any positive integers $m$ and $k$
$$
\A_n^{(m^k)}\subseteq \bigcup_{\{\lambda\vdash n|(m^k)\not \subseteq\lambda\}}
 C^\lambda .
$$

\smallskip

Note that Example 3.3 shows that the converse of Theorem 4.4 is false.

\smallskip

\noindent{\bf Proof.}
Let $\pi$ be a sequence of shape $\lambda$. By Greene's
Theorem, $\pi$ contains a disjoint union of $k$ increasing subsequences
of total size $\sum_{j=1}^k \lambda_j$. Denote this union by $\bar\pi$, and let
$\mu:=\sh(\bar \pi)$. Obviously, there are at most
$k$ parts in $\mu$ (i.e., $\mu=(\mu_1,\dots,\mu_k)$ with $\mu_k\ge 0$) and 
$\sum_{j=1}^k \mu_j=\sum_{j=1}^k \lambda_j$.
By Greene's Theorem,
$$
\sum_{j=1}^{k-1} \mu_j= \hbox{maximal size of a union of $k-1$ increasing subsequences in $\bar \pi$}\le
$$
$$
\le \hbox{maximal size of a union of $k-1$ increasing subsequences in $\pi$}=\sum_{j=1}^{k-1} \lambda_j.
$$
Hence, $\mu_k\ge \lambda_k$.
By assumption $(m^k)\subseteq \lambda$, so that
$m\le\lambda_k$. We conclude that there are exactly $k$ parts in $\mu$, and
$\mu_1\ge\cdots\ge\mu_k\ge m$.
In other words, $\mu_1'=k$ and $(k^m)\subseteq\mu'$.

Now, by the second part of Greene's Theorem,
$\bar \pi$ contains a disjoint union of $m$ decreasing subsequences
 of total size $km$.
 Denote this union by $\hat \pi$, and denote its shape by $\nu$.
$\hat\pi$ is a subsequence of $\bar\pi$, hence,
$$
\nu'_1=\hbox{length of maximal decreasing subsequence in $\hat\pi$} \le 
$$
$$
\le \hbox{length of maximal decreasing subsequence in $\bar\pi$}=
\mu'_1=k.
$$
On the other hand,
$$
|\nu|=\nu'_1+\dots+\nu'_m=km.
$$
This shows that the shape of the subsequence $\hat\pi$ is $\nu=(m^k)$.

\qed

\medskip

\section{General Shapes}

\begin{note}
\noindent{\bf Theorem 5.1.} For any  partition $\mu=(\mu_1,\dots,\mu_k)$ of $m$
and any positive integer $n$,
$$
\A_n^{\mu}\subseteq 
\bigcup_{\{\lambda\vdash n| (\mu_1^k)\not\subseteq \lambda\}} C^\lambda .
\leqno(5.1)
$$
\end{note}

\noindent{\bf Proof.} Let $\lambda$ be a shape such that 
$(\mu_1^k)\subseteq \lambda$.
By Theorem 4.4, any permutation of shape $\lambda$
contains a subsequence
of shape $(\mu_1^k)$. By Theorem 4.1, this subsequence contains a subsequence 
of shape $\mu$.
\qed

\medskip

Let $\av^\mu_n$ be the size of the set $\A_n^\mu$. 
Theorem 5.1 implies the following asymptotic estimates.

\begin{note}
\noindent{\bf Corollary 5.2.} 
For any fixed partition $\mu=(\mu_1,\dots,\mu_k)$, 
$$
\limsup_{n\rightarrow\infty}{(\av_n^\mu)^{1/2n}} \le
 \h(\mu)+\w(\mu) \leqno(5.2)
$$
and
$$
\max\{\h(\mu),\w(\mu)\}\le \liminf_{n\rightarrow\infty}{(\av_n^\mu)^{1/2n}} ,
\leqno(5.3)
$$
where the height of $\mu$ $\h(\mu):=\mu_1'-1$, and the width of $\mu$ 
$\w(\mu):=\mu_1-1$.
\end{note}

\noindent{\bf Proof.}
Let $\lambda$ be a partition of $n$, and let
$f^\lambda$ be the number of standard Young tableaux of shape $\lambda$.
By the Robinson-Schensted correspondence 
$$
(f^\lambda)^2=\#\{\pi\in S_n|\sh(\pi)=\lambda\}.
$$
Combining this fact with Theorem 5.1 we obtain
$$
\av_n^\mu\le \#\{\pi\in S_n|(\mu_1^k)\not\subseteq \sh(\pi)\}=
 \sum\limits_{\lambda\vdash n\wedge (\mu_1^k)\not\subseteq \lambda} 
 (f^\lambda)^2.
$$
The asymptotics of the sum on the right hand side was studied by Berele and
Regev [BR, Section 7]. By [BR, Theorem 7.21], for fixed $\mu_1$ and $k$
$$
 \sum\limits_{\lambda\vdash n\wedge (\mu_1^k)\not\subseteq \lambda} (f^\lambda)^2\sim c_1(\mu_1,k)\cdot n^{c_2(\mu_1,k)}\cdot
(\mu_1+k-2)^{2n}, \leqno(5.4)
$$
when $n$ tends to infinity. Here
 $c_1(\mu_1,k)$ and $c_2(\mu_1,k)$ are independent
of $n$. This proves the upper bound (5.2).

\smallskip

For the lower bound, note that by Schensted's Theorem any permutation
avoiding $(\mu_1)$ also avoids $\mu$.
Similarly, any permutation
avoiding $(1^k)$ also avoids $\mu$.
Thus
$$
\A_n^{(\mu_1)}\cup\A_n^{(1^{\mu_1'})}\subseteq\A_n^\mu.
$$
This implies that (for $n$ large enough; e.g., $n>(\mu_1-1)(\mu'_1-1)$ )
$$
\av_n^{(\mu_1)}+\av_n^{(1^{\mu_1'})}\le \av_n^\mu.
$$
Combining this inequality with (5.4) proves the lower bound (5.3). 
 
\qed

\medskip

\noindent {\bf Note:} For an evaluation of $\av_n^{(m)}$ for $m\le 4$ 
see [St Exer. 7.16(e)].
An asymptotic evaluation of
$\av_n^{(m)}$ for fixed $m>4$ was first done in [Re].

\medskip

\section{Other Special Cases}

\subsection{\bf Hooks}

In this subsection we show that for hook avoiding permutations
and $n$ large enough
the ``False Conjecture" is correct.

\smallskip

\begin{note}
\noindent{\bf Theorem 6.1.} For any hook $\mu=(m,1^{k-1})$ and  
$n> (2m-4)(2k-4)$
$$
\A^{(m,1^{k-1})}_n=\bigcup_{\{\lambda\vdash n|(m,1^{k-1})
 \not\subseteq \lambda\}} C^\lambda.
$$
\end{note}

\medskip

\noindent {\bf Note:} 
If either $m\le 3$ or $k\le 3$ then equality holds for all values of $n$.

\medskip

The following analogue of Lemma 4.3 characterizes permutations of hook shape.

\begin{note}
\noindent{\bf Lemma 6.2.}
A permutation $\pi$ has shape $(m,1^{k-1})$
if and only if $\pi$ is a union of an increasing subsequence of length $m$
and a decreasing subsequence of length $k$, intersecting in a unique element.
\end{note}

\noindent{\bf Proof.} By Schensted's Theorem, a permutation
 $\pi$ of shape $(m,1^{k-1})$
contains an increasing subsequence $\alpha$ with $|\alpha|=m$
and a decreasing subsequence $\beta$ with $|\beta|=k$, where
$|\alpha\cup\beta|\le |\pi|=m+k-1$.
Since necessarily $|\alpha\cap\beta|\le 1$, it follows that
 $|\alpha\cap\beta|=1$.

The converse follows similarly from Schensted's Theorem.

\qed

\medskip

\begin{note}
\noindent{\bf Lemma 6.3.}
Let $m$ and $k$ be positive integers.
\begin{itemize}
\item[(a)] If either $m\le 3$ or $k\le 3$ then every permutation whose shape contains the hook $(m,1^{k-1})$ has a subsequence of shape $(m,1^{k-1})$.
\item[(b)] If $m\ge 4$ and $k\ge 4$ then every permutation whose shape contains the hook 
$(2m-3,1^{k-1})$ or the hook $(m,1^{2k-4})$ has a subsequence of shape $(m,1^{k-1})$.
\item[(c)]  
For any $m\ge 4$ and $k\ge 4$ there exists a permutation whose shape 
contains $(2m-4,1^{2k-5})$, but it
has no subsequence of shape $(m,1^{k-1})$.
\end{itemize}
\end{note}

\smallskip

\noindent{\bf Note:} The results in (a) and (b) above are best possible,
as far as the assumed size of a hook contained in the shape is concerned.
For (a) this is clear, and for (b) this is the content of (c).

\medskip

\noindent{\bf Proof.}
We shall prove (b); the proof of (a) is similar.

\noindent{\bf (b)} Let $\pi$ be a permutation whose shape contains the hook
$(2m-3,1^{k-1})$, with $m,k\ge 4$.
Then $\pi$ has an increasing subsequence $\alpha$ of length $2m-3$
and a decreasing subsequence $\beta$ of length $k$. 
If $\alpha$ and $\beta$ intersect
(necessarily in a unique element), then by truncating $\alpha$ to $m$
elements we get by Lemma 6.2 a subsequence of shape $(m,1^{k-1})$.
Otherwise (i.e., assuming that $\alpha$ and $\beta$ do not intersect)
we will show that the union of $\alpha$ and $\beta$ contains the required subsequence.

Let $\alpha=(\alpha_1,\dots,\alpha_{2m-3})$ and $\beta=(\beta_1,\dots,\beta_k)$, so that
$\alpha_1<\dots<\alpha_{2m-3}$ and $\beta_1>\dots>\beta_k$.

Let ${\sl ind}(\alpha_i)$ denote the index of $\alpha_i$
in the union of $\alpha$ and $\beta$ (as a subsequence of $\pi$);
similarly for ${\sl ind}(\beta_j)$.

Concerning the element $\alpha_{m-1}$ there are three possibilities:
\begin{itemize}
\item[(1)] There is an index $1\le j\le k-1$ such that
$$
{\sl ind}(\beta_j)<{\sl ind}(\alpha_{m-1})<{\sl ind}(\beta_{j+1}).
$$
\item[(2)] ${\sl ind}(\alpha_{m-1})<{\sl ind}(\beta_1).$
\item[(3)] ${\sl ind}(\alpha_{m-1})>{\sl ind}(\beta_k).$
\end{itemize}
We shall deal with case (1); the other cases are similar.
Since $\beta_j>\beta_{j+1}$, there are now three subcases:
\begin{itemize}
\item[(1a)] $\beta_j>\alpha_{m-1}>\beta_{j+1}.$
\item[(1b)] $\alpha_{m-1}<\beta_{j+1}.$
\item[(1c)] $\alpha_{m-1}>\beta_j.$
\end{itemize}
In case (1a), $\alpha_{m-1}$ may be added to the decreasing subsequence $\beta$, to obtain two intersecting monotone subsequences of lengths $2m-3$ 
and $k+1$. By truncating these subsequences we will get an
increasing subsequence of length $m$ intersecting a decreasing subsequence of length $k$.

In case (1b), $(\alpha_1,\dots,\alpha_{m-1},\beta_{j+1})$
is an increasing subsequence of length $m$ intersecting $\beta$.

In case (1c), $(\beta_j,\alpha_{m-1},\alpha_m,\dots,\alpha_{2m-3})$
is an increasing subsequence of length $m$ intersecting $\beta$.

By Lemma 6.2, in all cases we obtain a subsequence of $\pi$ having shape $(m,1^{k-1})$.

\smallskip

\noindent{\bf (c)} The construction extends Example 3.2 (for which $m=k=4$):
take $\pi=(\gamma,\alpha,\delta,\beta)$, where
$\alpha$ and $\delta$ are increasing sequences of length $m-2$ and $\beta,
\gamma$ are decreasing sequences of length $k-2$:
$$
\alpha=(1,\dots,m-2); \qquad \beta=(m+k-4,\dots,m-1);
$$
$$ \gamma=(m+2k-6,\dots,m+k-3); \qquad
\delta=(m+2k-5,\dots,2m+2k-8).
$$
It is easy to see that an increasing subsequence of $\pi$ intersecting $\gamma$
must be contained (omitting
the intersection element itself) in $\delta$, so that its total length
is at most $m-1$.
Similar analysis of $\beta$ 
shows that an increasing
subsequence of length $m$ in $\pi$ must be contained in $(\alpha,\delta)$.
Analogously, a 
decreasing subsequence of length $k$ must be contained in $(\gamma,\beta)$.
The two subsequences cannot intersect.

\qed

\medskip

\noindent{\bf Proof of Theorem 6.1.} 
By Schensted's Theorem,
if a permutation $\pi$ has a subsequence of shape $(m,1^{k-1})$ then
it has an increasing subsequence of length $m$ and a decreasing subsequence
of length $k$. On the other hand, 
a permutation in 
$\bigcup_{\{\lambda\vdash n|(m,1^{k-1})\not\subseteq \lambda\}} C^\lambda$
has either no increasing subsequence of length $m$ or 
no decreasing subsequence of length $k$. 
Thus,
$$
\bigcup_{\{\lambda\vdash n|(m,1^{k-1})\not\subseteq \lambda\}} C^\lambda
\subseteq \A^{(m,1^{k-1})}_n.
$$

For the other direction,
assume that $\pi\in C^\lambda$ with $(m,1^{k-1})\subseteq \lambda$.
 Hence, $\lambda_1\ge m$ and $\lambda_1'\ge k$.
If either $m\le 3$ or $k\le 3$ then, by Lemma 6.3(a), $\pi$ has a subsequence of shape 
$(m,1^{k-1})$. Otherwise (i.e., if $m\ge 4$ and $k\ge 4$),
by assumption $(2m-4)(2k-4)<n=|\lambda|\le \lambda_1\cdot \lambda'_1$,
and therefore either $\lambda_1>2m-4$ or $\lambda'_1>2k-4$.
We can now use Lemma 6.3(b).
 
\qed

\medskip

\begin{note}
\noindent{\bf Corollary 6.4.} For any pair of positive integers
$m$ and $k$, and for $n\ge 4mk$
$$
\av^{(m,1^{k-1})}_n=\av^{(m)}_n+\av^{(1^k)}_n=
\sum\limits_{\lambda\vdash n\wedge\lambda_1<m} (f^\lambda)^2+
\sum\limits_{\lambda\vdash n\wedge\lambda'_1<k} (f^\lambda)^2,
$$
where $f^\lambda$ is the number of standard Young tableaux of shape $\lambda$.
\end{note}

Combining Corollary 6.4 with (5.4) we obtain
\begin{note}
\noindent{\bf Corollary 6.5.}
$$
\lim_{n\rightarrow\infty} (\av^{(m,1^{k-1})}_n)^{1/2n} = \max\{m-1,k-1\}.
$$
\end{note}

\medskip

\subsection{Avoiding $(2^2)$}

In this subsection we compute $\av^{(2^2)}_n$
and show that 
$$
\lim_{n\rightarrow \infty}(\av^{(2^2)}_n)^{1/2n} = \sqrt{2+\sqrt 2}.
$$
In particular, unlike the case of hooks, neither the lower bound nor 
the upper bound of Corollary 5.2 gives the correct limit in this case.

\medskip

\noindent
Example 3.3 shows that for any $n\ge 5$, 
$$
 \bigcup_{\{\lambda\vdash n|(2^2)
\not\subseteq\lambda\}} C^\lambda \not\subseteq\A_n^{(2^2)}.
$$
However, the opposite inclusion does hold.

\begin{note}
\noindent{\bf Proposition 6.6.} For any positive $n$, 
$$
\A_n^{(2^2)}\subseteq \bigcup_{\{\lambda\vdash n|(2^2)
\not\subseteq\lambda\}} C^\lambda .
$$
\end{note}

\smallskip

Proposition 6.6 is a special case of Theorem 4.4. Here we suggest an
independent and more informative proof of this result.

\smallskip

\noindent{\bf Proof.} By induction on $n$.
The claim obviously holds for $n\le 4$.
Assume that it holds for $n-1$, for some $n\ge 5$.

For the induction step observe that $C^{(2^2)}=\{2143,2413,3142,3412\}$ consists 
of all permutations in $S_4$
for which 1 and 4 are in the `middle'.
It follows that for any permutation $\pi$ in $S_n$,
if $\pi_1\not\in\{1,n\}$ and $\pi_n\not\in\{1,n\}$ then $\pi$
is not $(2^2)$-avoiding.
Therefore, if $\pi\in S_n$ is $(2^2)$-avoiding then either

$\pi_1\in \{1,n\}$ or $\pi_n\in \{1,n\}$.
Assume that $\pi_1\in \{1,n\}$. By the induction hypothesis
the shape of the subsequence $(\pi_2,\dots,\pi_n)$ 
does not contain $(2^2)$ and is therefore a hook
$(r,1^{n-r-1})$ for some $1\le r\le n-1$.
Adding $\pi_1=1$ increases the size of the
longest increasing subsequence by 1;
thus, by Schensted's Theorem the resulting shape is $(r+1,1^{n-r-1})$.
Adding $\pi_1=n$ increases the size of the 
longest decreasing subsequence by 1;
again, by Schensted's Theorem the resulting shape is $(r,1^{n-r})$.
The case $\pi_n\in \{1,n\}$ is similar.

\qed

\smallskip

\begin{note}
\noindent{\bf Corollary 6.7.} For any positive integer $n$
$$
\av_n^{(2^2)}={1\over 2}(2+\sqrt 2)^{n-1}+{1\over 2}
(2-\sqrt 2)^{n-1}.
$$
\end{note}

\noindent{\bf Proof.}
It follows from the proof of Proposition 6.6 that 
$$
\av_n^{(2^2)}=
4\cdot\av_{n-1}^{(2^2)}-2\cdot\av^{(2^2)}_{n-2}.
$$
The solution of this linear recursion 
(with appropriate initial values) gives the desired result.

\qed

\bigskip
\bigskip

\section{Final Remarks and Open Problems}

\subsection{Algebraic Structure}
Let $R$ be the set of all representatives of minimal length of left cosets 
of $S_m$ in $S_n$ (length here, as usual, is in terms of the Coxeter 
generators, i.e., adjacent transpositions).
For any partition $\mu$ of $m$, the set $C^\mu$ of all permutations of 
shape $\mu$ is a two-sided Kazhdan-Lusztig cell in $S_m$. For any $n\ge m$ 
the set of all permutations in $S_n$ which are not $\mu$-avoiding
coincides with the set $RC^\mu R^{-1}$. Theorem 5.1 claims that for
hook shapes the set $RC^\mu R^{-1}$ is a union of two-sided Kazhdan-Lusztig  
cells.
This phenomenon generalizes a beautiful well-known fact:
The set $RC^\mu$ (or: $C^\mu R^{-1}$) is a union of 
Kazhdan-Lusztig left (resp. right) cells 
[Sr, BV Prop. 3.15]. See also [GaR, Ro]. Barbasch and Vogan gave an algebraic 
proof of this fact by associating the set $RC^\mu$  to induced representations.
An algebraic interpretation for the results in this paper is required.
These and other relations with representation theory deserve further study.

\bigskip

\subsection{Asymptotics}
 Regev calculated, by considering Schensted's Theorem, the exact asymptotics of
$\av_n^{(m)}$ [Re].
In this paper we have generalized this ``RSK approach" to prove that
for any partition $\mu$ there exists a constant $c(\mu)$
such that, for any $n$,
$$
\av_n^\mu\le c(\mu)^n.
$$
Note that from Corollary 5.2 and Corollary 6.7 it also follows that, for 
$\mu$ not strictly contained in $(2^2)$, there exists a constant 
$\tilde c(\mu)>1$ such that $\av_n^\mu\ge\tilde c(\mu)^n$ for $n$ large enough.

A far reaching generalization was conjectured by Stanley and Wilf [Bo1].

\begin{note}
\noindent{\bf The Stanley-Wilf Conjecture.} 
For any fixed permutation $\sigma$ there exists a constant
$c(\sigma)$ such that, for any $n$
$$
\av_n(\sigma)\le c(\sigma)^n,
$$
where $\av_n(\sigma)$ is the number of all $\sigma$-avoiding
permutations in $S_n$.
\end{note}

\noindent
By a result of Arratia [Ar], if this conjecture holds then actually the limit
$\lim_{n\rightarrow \infty}\av_n(\sigma)^{1/n}$ always exists (and is finite).

The Stanley-Wilf conjecture holds for all $\sigma\in S_3$ [K, p. 238] and all
$\sigma\in S_4$ \cite{B1, B2}, as well as for many other cases (see [SSi], 
[Bo3] and their references). Recently, Alon and Friedgut [AF] have applied 
Davenport-Schinzel sequences to prove a somewhat weaker version of the 
conjecture for arbitrary $\sigma$.  An interesting challenge is to apply the 
``RSK approach" to attack the Stanley-Wilf Conjecture; namely,
to apply Greene's Theorem and methods presented in this paper
to sets avoiding a single permutation.

\bigskip

\noindent{\bf Acknowledgments.} The authors thank
Noga Alon, Miklos B\'ona, Ehud Friedgut, Nati Linial, Alek Vainshtein
and Julian West for useful discussions.
Special thanks to Amitai Regev for stimulating comments.

\end{document}